# The Dynamics of Deforming Manifold: A Mathematical Model

Xiaodong Zhuang[1], Nikos E. Mastorakis[2]
[1]Electronic Information College, Qingdao University, 266071 China
[2]Technical University of Sofia, Industrial Engineering Department, Kliment Ohridski 8, Sofia, 1000 Bulgaria

*Abstract*—In order to meet the requirements of practical applications, a model of deforming manifold in the embedded space is proposed. The deforming vector and deforming field are presented to precisely describe the deforming process, which have clear physical meanings. The proposed model is a modification of the general differential dynamic model, with constraints of spatial and temporal continuity on the deforming field. The deformation integral and derivative are presented as compact expressions of manifold deforming process. Moreover, a specific autonomous deforming field with flattening effect is defined, which provides a novel geometric viewpoint on data dimension reduction. The effectiveness of this autonomous deforming field is proved by numerical computation simulations, which indicate the promising potential of the proposed model in practical dimension reduction tasks.

## I. Introduction

MANIFOLD is a fundamental concept in modern mathematics. It plays an important role in modern researches in science and technology [1-20]. Especially, in the field of data science and artificial intelligence, manifold-based methods have been extensively studied [21-31]. There have been many researches on the mathematical description and topology property of manifold [32,33]. In mathematics, research mainly concentrates on static manifold. For deformed manifold, topological invariants provides global description under some topology-preserving deformation restrictions [33-35].

On the other hand, in practical applications, there are many cases of time-varying manifold with deforming shapes, such as a part of undulating sea surface, the time-varying twisted space caused by the gravity of a moving nearby celestial body, or the manifold of a continuous speech signal in the signal space. The study on such non-static manifold needs a quantitative and meticulous mathematical description of the manifold deforming process.

In this paper, the idea of deferential dynamics is introduced into the quantitative description of manifold deforming process. The deforming field on the manifold is proposed in the embedded space, which derives integral and deferential descriptions of the deforming process. A specific autonomous deforming field is also proposed to flatten the curved manifold in the embedded space, which has potential application in dimension reduction in data analysis. Simulation experiments prove the effectiveness of the proposed deforming field. The proposed model gives new insights into the geometric aspects of data dimension reduction.

## II. The Deforming Field on Manifold and the Deformation Integral

In mathematics, the theory of deferential dynamics accurately describes the evolution of a system over time [36,37]. The trajectory of a point in the phase space can be determined by a group of deferential equations. The deforming process of a manifold can be regarded as the deformation of the corresponding geometry in the embedded space. In the embedded space, only the points on the manifold are involved in the deformation process. However, in differential dynamics, all the points in the phase space are involved in the "stream" determined by the solutions of the deferential equations. Moreover, there are infinite points on a smooth manifold, and its deformation involves all its points moving synchronously. This is different from the general deferential dynamic system, where each point independently forms its own trajectory in the phase space. Therefore, it is necessary to define specific mathematical descriptions for manifold deformation process, i.e. the dynamics of deforming manifold.

There are two different forms of describing the manifold mathematically: the intrinsic form and extrinsic form [38,39].

The intrinsic form does not rely on the specific coordinate system of the embedded space, and the corresponding non-Euclidean geometry plays a fundamental role in the description of time-space in modern physics. On the other hand, data vectors in modern data science have a natural form of extrinsic representation, which are the coordinates of data manifold points in the embedded space. Therefore, it is natural to represent the data manifold deformation by the corresponding deforming geometry in the embedded space. In this paper, the extrinsic definition of manifold deforming dynamics is proposed.

In practical data mining and machine learning tasks, the data manifold is usually represented by a finite data set, which consists of discrete samples from the data manifold. Therefore, the dynamics of deforming manifold is defined for bounded manifold with borders (one simple example is a pack of sea surface undulating randomly). Let $M$ denote the deforming manifold. $p$ is a point on $M$. In the deforming process, the moving of $p$ is expressed by the change of its coordinates in the embedded space. The embedding of manifold $M$ into $R^n$ is homeomorphism. If $M$ is bounded with borders, it corresponds to a bounded geometry $M$ in $R^n$. Suppose $M$ continuously deforms in the time interval [0, T]. The corresponding geometry also deforms synchronously in $R^n$.

In the deforming process, $M$ changes its shape in $R^n$, and each point $p$ on $M$ also changes its position in $R^n$ simultaneously. It should be noted that $p$ is a geometric element on $M$, while its position (coordinates in $R^n$) is a time-variant function $z_p(t)$. In the model proposed, for each point $p$ on $M$, there is a vector-valued function $\vec{v}(p,t)$. At any time $t_0$, all these vectors constitute a vector field $D_M(t_0)$ on $M$. If $\vec{v}(p,t)$ and $D_M(t_0)$ satisfy:

1. For any point $p_0$ on $M$, $\vec{v}(p_0,t)$ is a continuous vector function of time $t$ (the constraint of temporal continuity).
2. For any time $t_0$, $D_M(t_0)$ is continuous on $M$ (the constraint of spatial continuity on $M$).

Then $D_M(t)$ is defined as a continuous deforming field, which is a time-varying vector field on $M$. And $\vec{v}(p,t)$ is defined as the deforming vector. $D_M(t)$ can be regarded as a specific vector bundle defined on $M$, where there is a one-dimensional vector space with $\vec{v}(p,t)$ as the base vector on each point $p$. In the proposed model, a deforming field $D_M(t)$ can determine the trajectory of each point $p$ in a deforming process:

$$\frac{d\vec{z}_p}{dt} = \vec{v}(p,t) \qquad (1)$$

where $z_p$ denotes the coordinates of $p$ in $R^n$. $\vec{v}(p,t)$ is the value of $D_M(t)$ on point $p$. $\vec{v}(p,t)$ can be physically explained as the velocity vector of $p$ in the deforming process. Although Equation (1) has the general form of differential equations of the first order standard type, $\vec{v}(p,t)$ has the constraints of temporal and spatial continuity, which are necessary in a continuous deforming process. And the constraint of spatial continuity on $M$ is the key difference between the proposed model and the general differential dynamic systems.

For any point $p$ on $M$, the solution to Equation (1) gives the trajectory of $p$ in $R^n$ in the deforming process. Moreover, at time $t_0$, all the points form a geometry representing the current "shape" of the manifold. Solving Equation (1) for all the points on $M$ is defined as the "deformation integral", which is expressed in a compact form:

$$M(t) = \int_o^t D_M(\tau) \cdot d\tau \qquad (2)$$

where $D_M(t)$ is the deforming field, and $M(t)$ is the manifold geometry in embedded space at time $t$. Equation (2) is a compact representation of the deforming process under the deforming field $D_M(t)$. Equation (1) is from the viewpoint of a single point on the manifold, while Equation (2) is from the viewpoint of the whole manifold. Because the deformation of the whole manifold is eventually composed of all the points' movement, Equation (1) explains the actual meaning of Equation (2). Because the derivation and integral are inverse operation, Equation (2) can be written in another form:

$$D_M(t) = \frac{dM(t)}{dt} \qquad (3)$$

which is defined as the deformation derivative.

## III. AUTONOMOUS DEFORMATION AND TOPOLOGY-PRESERVING DEFORMATION

### A. The Autonomous Deformation

In analogy with the autonomous system in deferential dynamic system, the autonomous deformation is defined here. In a deforming field $D_M(t)$, if the deforming vector $\vec{v}(p,t)$ at point $p$ on $M$ is only determined by the current "shape" of $M$ and the relative position of $p$ on $M$, the deforming process is defined as the autonomous deformation. From a global viewpoint, the deforming field at time $t$ is only determined by the manifold shape at that time:

$$D_M(t) = F(M(t)) \qquad (4)$$

Equation (4) is from the viewpoint of the whole manifold. $F$ is a mapping (or functional) that represents the determination of $D_M(t)$ by $M(t)$. Substitute (4) into (2), the integral form of autonomous deformation is obtained as:

$$M(t) = \int_o^t F(M(\tau)) \cdot d\tau \qquad (5)$$

Similarly, substitute Equation (4) into (3), the differential form of autonomous deformation is obtained:

$$\frac{dM(t)}{dt} = F(M(t)) \qquad (6)$$

Equation (6) implies that, how the manifold will deform at time $t$ (i.e. the deforming field $D_M(t)$) only depends on the manifold shape at time $t$. Equation (5) implies that, if $F$ is given, the autonomous deformation is only determined by the initial status of $M$ (i.e. the manifold initial shape). This is the unique property of autonomous deformation, i.e. the self-evolution of the manifold without external affection. In spite of the extrinsic definition of the autonomous deforming field, it is actually intrinsic given the mapping $F$, just like the intrinsic curvature of a surface with an extrinsic form of definition.

### B. The Topology-Preserving Deformation

Among all the possible deforming process of a manifold, the most practically important ones are those keeping the manifold's topological structure unchanged. This kind of deforming is defined as the topology-preserving deformation. If $M$ keeps its topology structure throughout the deforming process, two different points $p_1$ and $p_2$ on $M$ will not move to one position in $R^n$ at the same time (i.e. different points will not adhere to each other). And fragmentation will not happen in a continuous region on $M$.

Neighborhood is the basis of topological structure [34,35,40]. In the proposed model, a convenient way of neighborhood definition for points on the manifold is proposed here. The following definition is for the bounded smooth manifold mentioned in section II, which is a common situation in practical applications. In another word, the dataset only occupies limited area in the data space, and the "shape" of the manifold is not extensively complex (i.e. there is no infinitely dense folding on the manifold).

Manifold $M$ corresponds to a high dimensional geometry in $R^n$, and each point $p$ on $M$ has $n$-dimensional coordinate in $R^n$. In dimension reduction tasks, the initial shape of the data manifold $M$ is especially important, because it determines the topological structure. The neighborhood relationship between any two points on $M$ is judged according to the initial dataset. In the deforming process, if the neighborhood relationship between any two points keeps the same as that on the initial manifold, the topological structure of the manifold can be regarded as unchanged. Therefore, here the neighborhood of a point is determined according to the initial manifold's shape $M(0)$.

The $r$-neighborhood of point $p$ is defined as the set of points on $M(0)$ whose Euclidean distance to $p$ is less than $r$. The Euclidean distance is measured in the embedded space $R^n$, which is convenient to compute. For the convenience of later discussions, the $r$-neighborhood is defined as deleted neighborhood:

$$U_p^r \triangleq \{q | q, p \in M(0), 0 < dist(q,p) < r\},\ r > 0 \quad (7)$$

where $dist(q,p)$ is the Euclidean distance between $q$ and $p$. Because smooth manifold is locally homeomorphic in $R^n$, there exists $r$ with sufficiently small value that makes the above definition reasonable. On the initial shape $M(0)$, each point $p$ has an $r$-neighborhood $U_p^r$, which is a patch on $M(0)$. Preserving the neighborhood relationship is important for keeping the manifold's topological structure. It is a necessary condition for topology-preserving deformation that the neighborhood relationship for all the points on $M$ are preserved in the deforming process.

IV. AN AUTONOMOUS DEFORMING FIELD WITH FLATTENING EFFECT

A data vector can be regard as a sample from the data manifold. In data analysis applications, data vectors are usually of high dimension, such as digital images of high resolution. The element number of the data vector is usually not consistent with the intrinsic dimension of the data manifold. The data manifold usually has a much lower intrinsic dimension (imagine a piece of paper is wadded into spitball, which is an embedding into the 3D space). In many dimension reduction methods, it is necessary to estimate the intrinsic dimension of the data manifold.

From the viewpoint of manifold deformation, nonlinear dimension reduction is geometrically equivalent to the flattening of the manifold in the embedded space $R^n$. In such flattening process, the original neighborhood relationship between manifold points should be preserved, while the non-neighbor points should separate as far as possible (Imagine again the process of flattening the wadded paper). Inspired by this geometrical analogy, an autonomous deforming field is proposed, which preserves the original neighborhood relationship between data points, while separates non-neighborhood points far away. The intrinsic dimension of the data manifold can also be revealed at the same time.

The following definition is for the bounded manifold mentioned in section II. Suppose on the initial manifold $M(0)$, a point $p$ has a neighborhood $U_p^r$. $q$ is one of $p$'s neighbor point on $M(0)$ (i.e. $q \in U_p^r$), but $w$ is not $p$'s neighbor point. When $t>0$, $M$ deforms continuously. And the three points $p, q, w$ also move in the embedded space (i.e. there coordinates are functions of $t$: $p(t), q(t), w(t)$).

In order to flattening the manifold in $R^n$, an autonomous deforming field is defined, which separates non-neighbor points while maintain the distance between neighbor points as the original. To maintain the distance between neighbor points, an "elastic" interaction is defining between any two neighbor points $p$ and $q$ ($q \in U_p^r$):

$$K_1 \cdot \frac{\vec{p}(t) - \vec{q}(t)}{|\vec{p}(t) - \vec{q}(t)|} \cdot (|\vec{p}(0) - \vec{q}(0)| - |\vec{p}(t) - \vec{q}(t)|) \quad (8)$$

where $K_1$ is a constant and $K_1>0$. $\vec{p}(t)$ and $\vec{q}(t)$ are the coordinates (i.e. position vectors) of the two neighbor points at time $t$. $|\vec{p}(0) - \vec{q}(0)|$ is the initial Euclidean distance between $p$ and $q$, while $|\vec{p}(t) - \vec{q}(t)|$ is the Euclidean distance between them at time $t$. By such definition, when the two neighbor points $p$ and $q$ get far away from each other, the vector of elastic interaction has the direction that lets $p$ get closer to $q$. Otherwise, this vector has the direction that separates $p$ from $q$. Therefore, the elastic interaction has the effect of keeping the original distance unchanged between neighbor points.

On the other hand, in order to make the non-neighbor points separate from each other, the repelling interaction is defined between two non-neighbor points $p$ and $w$:

$$K_2 \cdot \frac{\vec{p}(t) - \vec{w}(t)}{|\vec{p}(t) - \vec{w}(t)|} \quad (9)$$

where $K_2$ is a constant and $K_2>0$. $\vec{p}(t)$ and $\vec{w}(t)$ are the coordinates of the two non-neighbor points at time $t$. The vector of repelling interaction has the direction that always separates $p$ away from $w$.

The deforming vector in Equation (1) is then defined as the combination of the above two interactions:

$$\frac{d\vec{z}_p}{dt} = \vec{V}_e(t) + \vec{V}_r(t) \quad (10)$$

where

$\vec{V}_e(t) =$

$$\int_{U_p(t)} K_1 \cdot \frac{\vec{p}(t) - \vec{q}(t)}{|\vec{p}(t) - \vec{q}(t)|} \cdot (|\vec{p}(0) - \vec{q}(0)| - |\vec{p}(t) - \vec{q}(t)|) \cdot dq$$

(11)

and

$$\vec{V_r}(t) = \int_{M(t) - U_p(t)} K_2 \cdot \frac{\vec{p}(t) - \vec{w}(t)}{|\vec{p}(t) - \vec{w}(t)|} \cdot dw \qquad (12)$$

In Equation (11) and (12), $M(t)$ is the manifold at time $t$, and $U_p(t)$ is the area constituted by $p$'s neighbor points at time $t$. In the integral, $dq$ and $dw$ are the two infinitesimal volume elements on $M(t)$, which containing $q$ and $w$ respectively. In Equation (11) the integral obtains the overall effect of keeping the original distance unchanged between $p$ and its neighbor points. Similarly, in Equation (12) the integral obtains the overall effect of separating $p$ away from its non-neighbor points as far as possible.

## V. COMPUTER SIMULATION EXPERIMENTS

For data manifolds in practical applications, there are usually groups of data points, but no analytical expressions. Therefore, it is suitable to use numerical computation to study the proposed model of deforming manifold. Numerical computation has been successfully used in many research areas [41-43]. In this paper, numerical simulation is implemented to trace the autonomous deforming process under the deforming field of Equation (10).

The continuous manifold $M$ is spatially discretized into a mesh of sample points in $R^n$. The deforming vector $\vec{v}(p, t)$ is calculated discretely on the mesh, by approximating the integral in Equation (11) and (12) with discrete summation. Temporal discretization is also implemented for simulating the deforming process. For a small time interval $\Delta$, the manifold shape at $t+\Delta$ is determined by moving each mesh point according to the deforming vector on it. In the simulation, the two constants $K1$ and $K2$ in Equation (11) and (12) are 0.1 and 0.0002 respectively, with which stable and convergent simulation is achieved.

Simulation experiments have been carried out on a group of typical manifolds. The experimental results for a half circle as one-dimensional manifold embedded in $R^2$ are shown in Fig. 1 and Fig. 2. The radius of the half circle is 69, and the manifold is discretized into 129 sample points in the simulation. On the initial manifold, the neighborhood radius is set to 3.36. Fig. 1 demonstrates the deforming field defined by Equation (10) on the initial manifold (the half circle). It can be seen that the deforming vectors have a tendency of stretching the curve. Some intermediate results of the deforming process have been recorded in the simulation, which are shown in Fig. 2 as a demonstration of the deformation process. The arrow in Fig. 2 indicates the deformation sequence. Fig. 2 indicates that the initial curved manifold is flattened to a straight line under the autonomous deforming field.

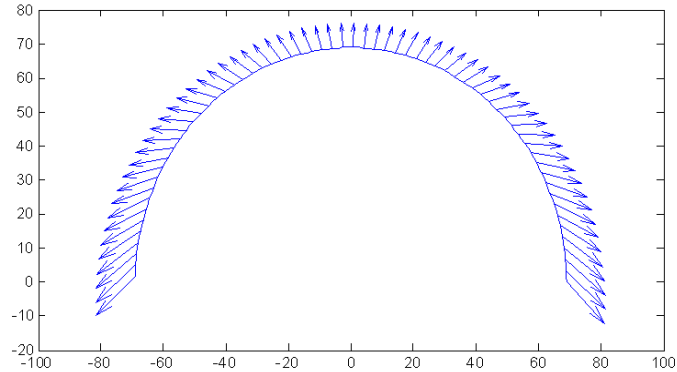

Fig. 1 The autonomous deforming field on the initial manifold of a half circle

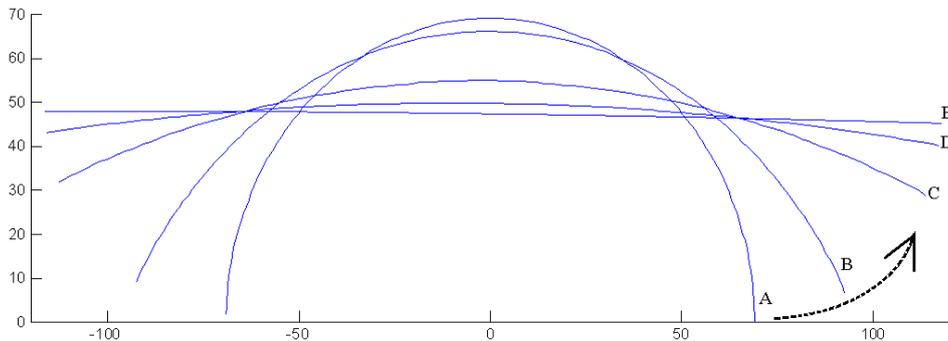

Fig. 2 The demonstration of deforming process of a half circle under the autonomous deforming field. (curve A is the initial manifold of a half circle; line E is the final deforming result; curve B, C, and D are three intermediate results in the deforming process)

Another simulation experiment is carried out for a manifold with more complex shape, a spiral line on the plane. The parametric equations of the spiral line is:
$$x(t) = (t + 10) \cdot cos(\pi t)$$
$$y(t) = (t + 10) \cdot sin(\pi t)$$
$$t \in [-1.0, 5.0]$$

The manifold is discretized into 600 sample points in the simulation. On the initial manifold, the neighborhood radius is set to 1.2. Fig. 3 demonstrates the deforming field defined by Equation (10) on the initial manifold. It can be seen that the deforming vectors have a tendency of expanding outwards. Intermediate results of deforming have been recorded in the simulation, which is shown in Fig. 4 as a demonstration of the deformation process. It can be seen that the spiral is gradually unfolded, and then gradually flattened under the autonomous deforming field of Equation (10).

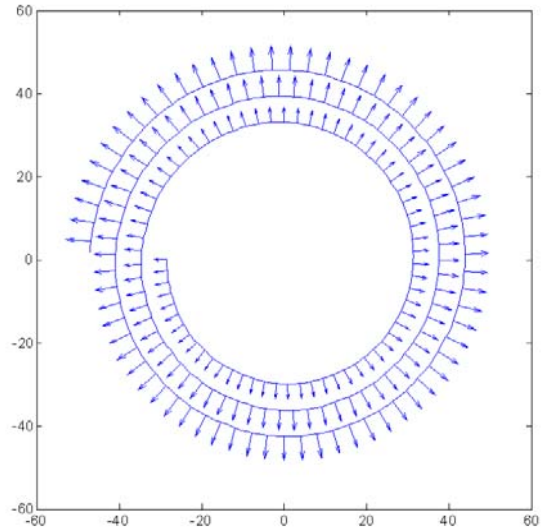

Fig. 3 The autonomous deforming field on the initial manifold of a spiral line

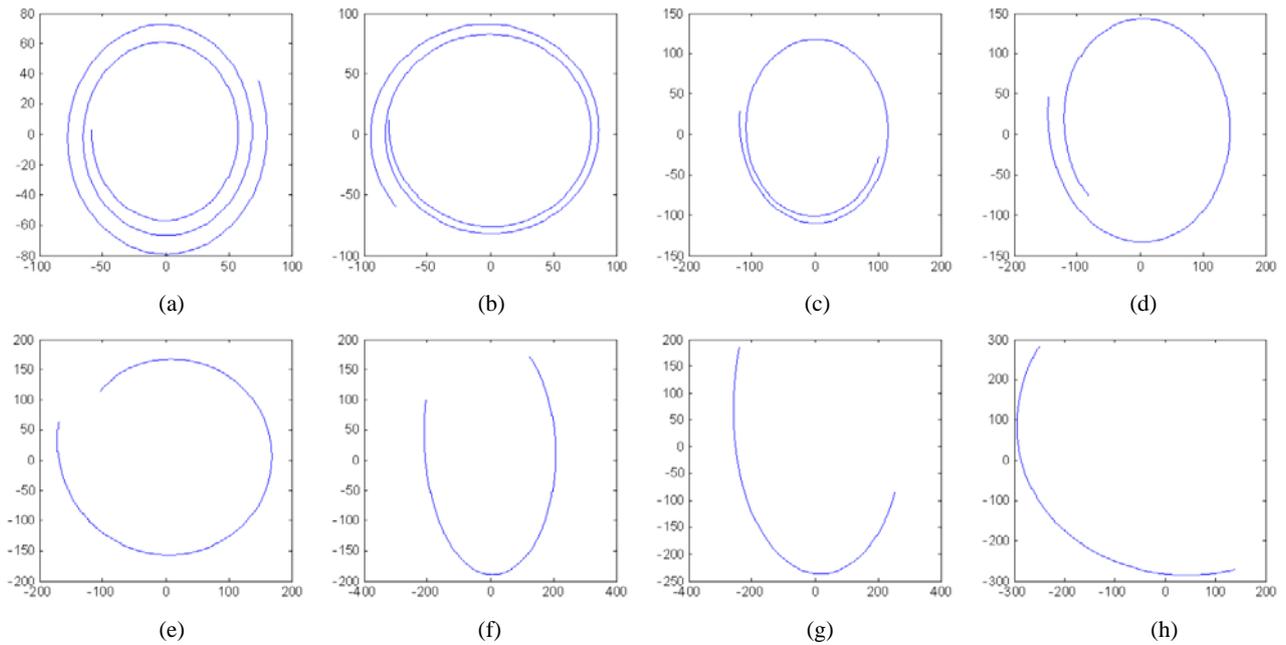

Fig. 4 The demonstration of deforming process of a spiral line under the autonomous deforming field by the intermediate result image sequence (a) to (h)

Simulation experiments have also been done for manifolds of higher dimension. The experimental results for the S-curve are shown in Fig. 5 to Fig. 8. The S-curve is a 2-D curved surface embedded in $R^3$, whose transverse section is a curve of the shape "S". The S-curve is discretized into 360 sample points in the simulation, and its discrete mesh is shown in Fig. 5(a). Fig. 6 demonstrates the deforming field defined by Equation (10) on the initial manifold. Fig. 6(b) shows the top view of the initial deforming field. It is clear that the vector directions of the deforming field have the effect of stretching the curved surface. Intermediate results have been recorded in the simulation of deforming, which are shown in Fig. 7 and Fig. 8 as a demonstration of the deformation process. These intermediate results are shown in two different angles of view: top view (in Fig. 8) and common 3D perspective (in Fig. 7). The results shows that higher dimensional manifold can also be flattened by the proposed autonomous deforming field. Because flattening implies dimension reduction, the proposed autonomous deformation mechanism has potential application in practical data analysis.

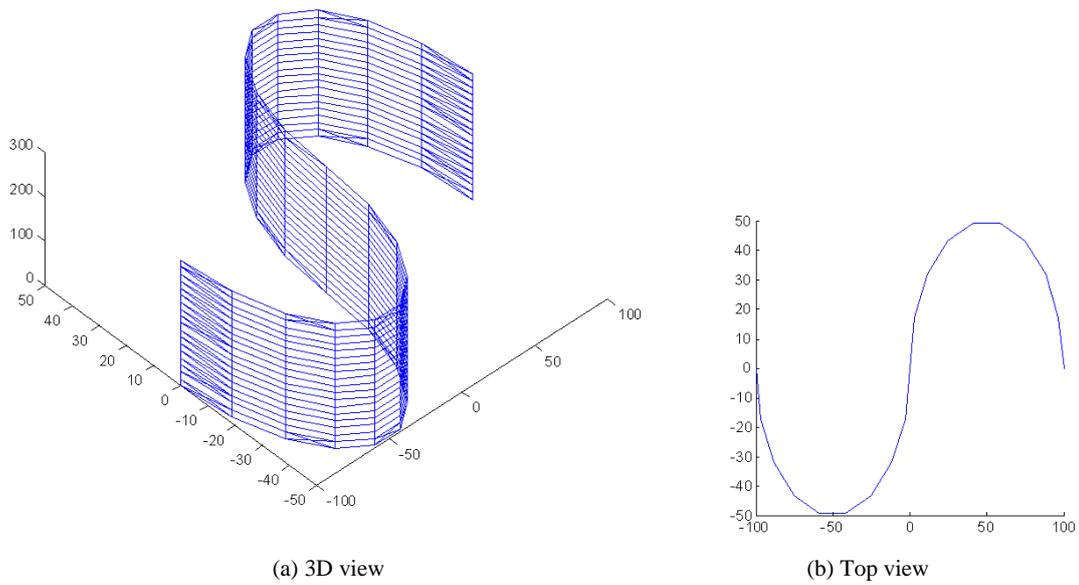

(a) 3D view  (b) Top view
Fig. 5 The initial manifold of S-curve

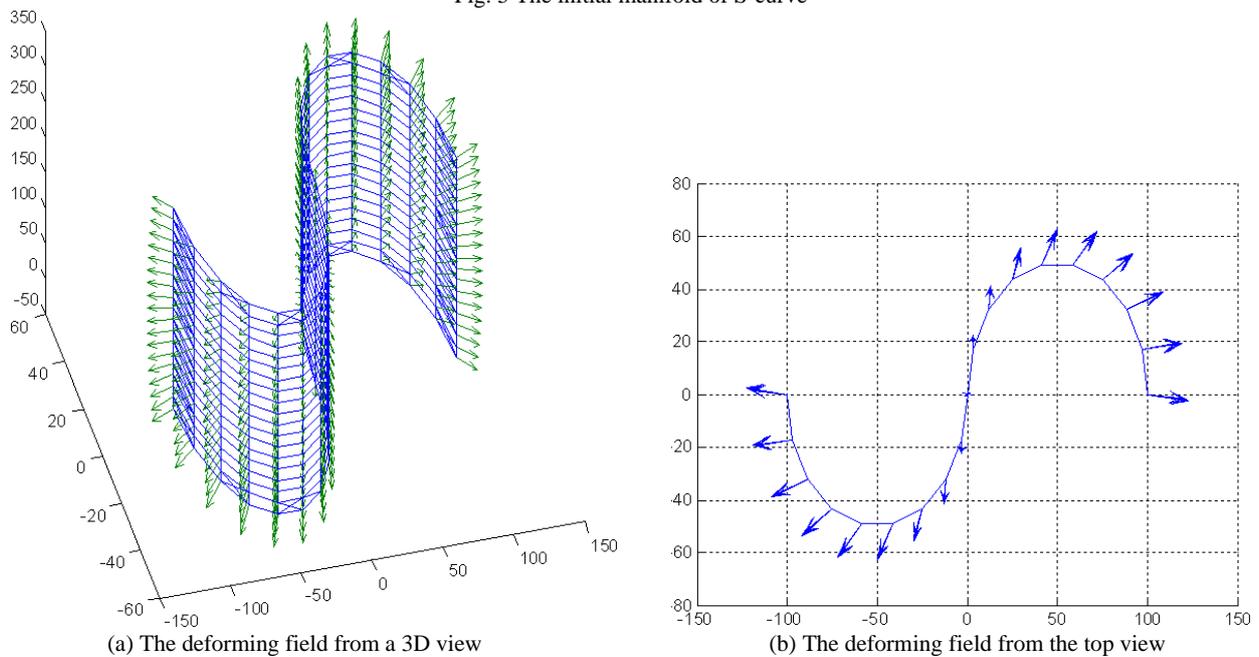

(a) The deforming field from a 3D view  (b) The deforming field from the top view

Fig. 6 The demonstration of the autonomous deforming field on the initial manifold of S-curve

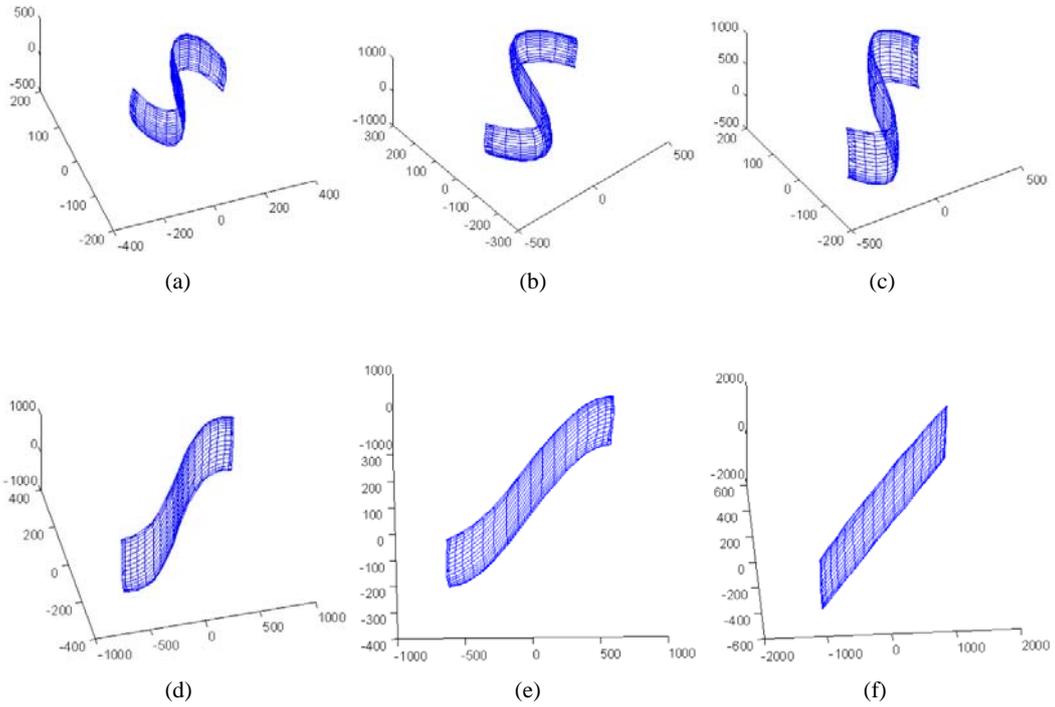

Fig. 7 The demonstration of S-curve deforming process (**3D view**) under the autonomous deforming field
(the intermediate result sequence (a) to (f))

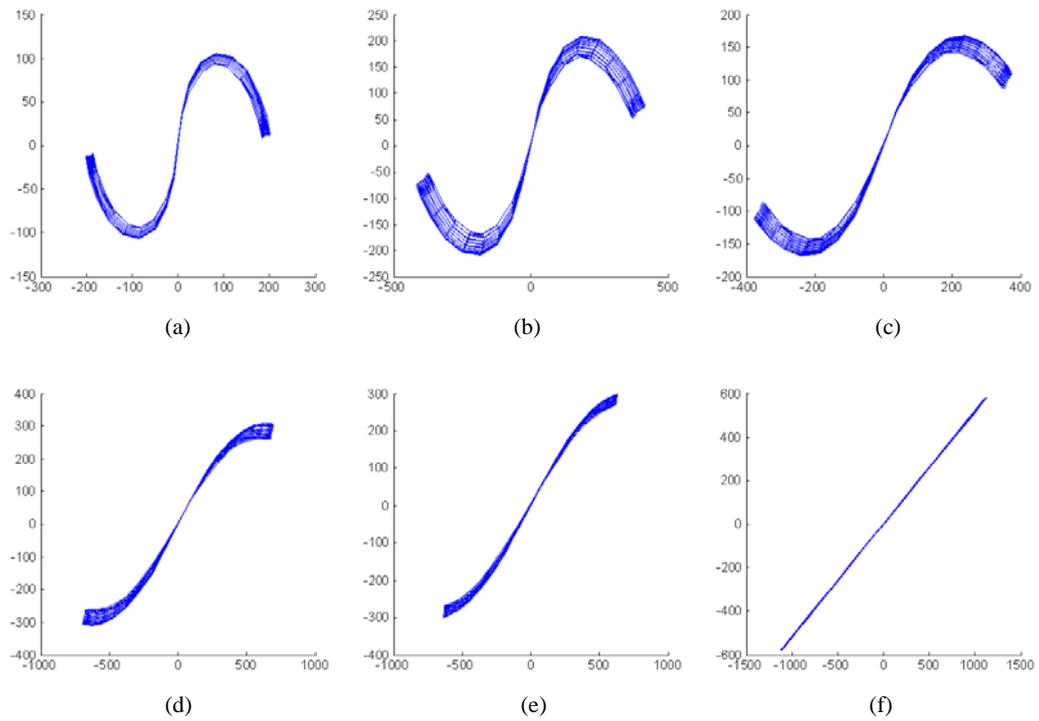

Fig. 8 The demonstration of S-curve deforming process (**top view**) under the autonomous deforming field
(the intermediate result sequence (a) to (f))

## VI. Conclusion

Traditional differential dynamics theory describes the movement of single points in the phase space, which can not directly represent the deformation of a whole manifold (i.e. the synchronously movement of a set of points under some topological constraints). To meet the practical requirement of quantitative description of manifold deformation, a model of deforming field is proposed, which implements the precise mathematical representation of the deforming process. Based on the deforming field, the differential and integral forms are presented to express the manifold deforming process, which provides a new starting point for further quantitative study of manifold deformation dynamics. Moreover, the autonomous deforming field is proposed as a model of manifold self-evolution. And a specific autonomous deforming field with flattening effect is defined, which is simulated by numerical computation on computer for some typical manifolds embedded in $R^2$ or $R^3$. The simulation results prove the effectiveness of the deforming field proposed, and indicate its potential in practical data analysis tasks. Future study will be carried out on mathematical properties of the deformation derivative and integral. The intrinsic property of the autonomous deforming field will be studied. Further application of the deforming field in practical data analysis (i.e. dimension reduction, feature extraction, etc.) will also be investigated.